\documentclass[11pt,a4paper]{article}

\usepackage{xcolor}
\usepackage[margin=28mm]{geometry}
\usepackage{amsmath,amssymb,amsthm}
\usepackage[colorlinks=true,linkcolor=blue,citecolor=blue,urlcolor=blue]{hyperref}

\newtheorem{thm}{Theorem}
\newtheorem{lem}[thm]{Lemma}
\newtheorem{conj}[thm]{Conjecture}
\newtheorem{remark}[thm]{Remark}
\newcommand{\al}{\alpha}
\newcommand{\om}{\omega}
\allowdisplaybreaks[1]
\usepackage[shortlabels]{enumitem}

\title{Independence number, clique number, and maximum degree}
\author{Xinmin Hou
\thanks{School of Mathematical Sciences, University of Science and Technology of China, Hefei 230026, China.
Hefei National Laboratory, University of Science and Technology of China, Hefei 230088, Anhui, China. Email: xmhou@ustc.edu.cn.}
\and 
Xiangyang Wang\thanks{
School of Mathematical Sciences, University of Science and Technology of China, Hefei 230026, China.
Email: wangxiangyang@mail.ustc.edu.cn.}
}
\date{}

\begin{document}
\maketitle

\begin{abstract}
We prove that for every integer $r \ge 3$, every $K_r$-free graph $G$ on $n$ vertices with minimum degree $\delta(G) > (5r-11)n/(5r-5)$ satisfies $\alpha(G) \ge n/(r-1)$. This resolves a conjecture posed by Brandt in 2001 and establishes the sharp minimum-degree threshold for guaranteeing an independent set of size at least $n/(r - 1)$ in $K_r$-free graphs. 
\end{abstract}

\section{Introduction} 
Turán's theorem, a fundamental and widely known result, states that the \(K_r\)-free graphs on \(n\) vertices with the maximum number of edges are exactly the balanced complete \((r-1)\)-partite graphs \cite{Turan}. It is natural to ask what structural properties are retained when the graph is sufficiently dense, though not necessarily extremal. The celebrated Andrásfai--Erdős--Sós theorem \cite{AES} is a classical result concerning the chromatic number of $K_r$-free graphs with large minimum degree (for a short proof, see Brandt \cite{Br02}). It can be viewed as a minimum-degree analogue of the Erdős--Simonovits stability theorem \cite{Er01,Er02,Simonovits1968} for the structure of \(K_r\)-free graphs with nearly the maximum number of edges.

\begin{thm}[\cite{AES}]\label{AES}
    Let $r\ge3$ be an integer, and let $G$ be a $K_r$-free graph on
$n$ vertices. If
\[
\delta(G)>\frac{3r - 7}{3r - 4}\,n,
\]
then $\chi(G)\le r - 1$.
\end{thm}

Analogously,  Brandt \cite{Br01} posed the following conjecture,  replacing the chromatic number with the independence number.

\begin{conj}[\cite{Br01}]\label{conj:main}
    Let $r\ge 3$ be an integer, and let $G$ be a $K_r$-free graph on
$n$ vertices. If
\[
\delta(G)>\frac{5r-11}{5r-5}\,n,
\]
then $\al(G)\ge n/(r-1)$.
\end{conj}

For $r = 3$, the result can be deduced from Theorem \ref{AES}, which shows that $G$ is bipartite. The following structural result was proved by Goddard and Lyle \cite{GL}.

\begin{thm}[\cite{GL}]
    Let $G$ be  a maximal $K_r$-free graph on n vertices with $r\ge 4$ and \[\delta(G) > \frac{2r - 5}{2r - 3}n.\] Then, $G$ is the join of a triangle-free graph and a complete 
    $(r - 3)$-partite graph; in particular, $G = \Gamma \bigvee{(\bigvee_{i = 1}^{r - 3}\overline{K_{a_i}})}$ for some $a_i$, where $\Gamma$ is a triangle-free graph with minimum degree larger than $|\Gamma|/{3}$. 
\end{thm}
We may assume that $G$ is maximal $K_r$-free, since adding edges
cannot decrease its minimum degree or increase its independence number. If $\Gamma$ is bipartite, then $G$ can be partitioned into $r-1$ independent sets; consequently, one of them has size at least $n/(r-1)$. Otherwise, $\Gamma$ is non-bipartite. Then, by Theorem \ref{AES}, we have $\delta(\Gamma) \le \frac{2}{5}|\Gamma|$. Hence,
\[
\delta(G) - \sum_{i=1}^{r-3} a_i \le \frac{2}{5}\left(n - \sum_{i=1}^{r-3} a_i\right).
\]
It follows that $\sum_{i=1}^{r-3} a_i \ge ((r - 3)/(r - 1))n$ for $r \ge 4$. Thus, there also exists an independent set of size at least $n/(r-1)$. Therefore, a minimum degree greater than $((2r - 5)/{(2r - 3)})n$ ensures the existence of such an independent set. For $r = 4$, $(2r-5)/(2r-3) = (5r-11)/(5r-5)$. Thus, Conjecture \ref{conj:main} is also true for $r = 4$. Actually, Brandt \cite{Br01} proved the conjecture for $r=4$ and $r=5$. More recently, Lyle \cite{Lyle} settled the conjecture for the case $r=6$ via a novel technique focusing on  the behavior of large independent sets in $G$. For general minimum degree conditions with $r\ge 6$, Brandt \cite{Br01} proved that $((r^2-5)/(r^2+r-2))n$ is sufficient; later, this was improved to $((r-3)/(r-2)+1/((r-2)^2(r-1)))n$ by Lyle \cite{Lyle14}. 

For two vertex-disjoint graphs $G$ and $H$, the
\emph{join} $G\vee H$ is the graph obtained from
their disjoint union by adding all edges between
$V(G)$ and $V(H)$. For a positive integer $t$, the \emph{blow-up} $G[t]$ of a graph $G$ is obtained by replacing each vertex $v\in V(G)$ with an independent set $V_v$ of size $t$, and joining every vertex of $V_u$ to every vertex of $V_v$ if and only if $uv\in E(G)$. The graphs below show that $((5r-11)/(5r-5))n$ is the optimal minimum degree threshold: they satisfy $\delta(G) = ((5r-11)/(5r-5))n$ for $r$ odd and $\delta(G)/n \to (5r-11)/(5r-5)$ as $a \to \infty$ for $r$ even, while $\alpha(G) < n/(r-1)$.
\begin{itemize}
    \item For every odd $r\ge 3$, let $t=(r-1)/2$ and consider the graph $G=\bigvee_{i=1}^{t}C_5$.
    \item For every even $r\ge4$ and every positive integer $a$, let $t=(r-2)/2$ and consider the graph $G=\overline{K_{5a - 1}}\vee\bigl(\bigvee_{i=1}^{t}C_5[2a]\bigr)$. 
\end{itemize}

In this paper, we establish the following stronger result. The proof integrates the techniques of Lyle \cite{Lyle} with the independent transversal approach utilized by King \cite{King2011}.

\begin{thm}\label{thm:counting}
Let $a>0$ be a real number and $G$ a graph on $n$ vertices. If
\[
\om(G)<a\qquad\text{and}\qquad\Delta(G)<\frac{6a}{5}-1,
\]
then $\al(G) > n/a$.
\end{thm}
Equivalently, the theorem states that every graph $G$ satisfies
\[
\alpha(G)\ge
\frac{|V(G)|}
{\max\{\omega(G),\,5(\Delta(G)+1)/6\}}.
\]
Cranston and Rabern~\cite{Cranston2015}
proved this bound for vertex-transitive graphs. We can then deduce Conjecture \ref{conj:main} from Theorem \ref{thm:counting}.

\begin{proof}[Proof of Conjecture~\ref{conj:main} from Theorem~\ref{thm:counting}]
Suppose instead that $\al(G)<n/(r-1)$. Put $H=\overline G$ and
$a=n/(r-1)$. Then $\om(H)=\al(G)<a$ and
\[
\Delta(H)=n-1-\delta(G)<\frac{6n}{5(r-1)}-1=\frac{6a}{5}-1.
\]
Theorem~\ref{thm:counting} and $\al(H)=\om(G)\le r-1$ now imply
\[
n<a\al(H)\le a(r-1)=n,
\]
a contradiction.
\end{proof}

Reed~\cite{Reed1998} conjectured that every graph $G$ satisfies
$
\chi(G)\le
\left\lceil\frac{\Delta(G)+1+\omega(G)}{2}\right\rceil.
$
For the independence number, Fajtlowicz~\cite{Fajtlowicz1978}
proved that
$
\alpha(G)\ge \frac{2n}{\omega(G)+\Delta(G)+1}.
$
Under the hypotheses of Theorem~\ref{thm:counting}, we have
$\omega(G)+\Delta(G)+1<11a/10$, so Fajtlowicz's bound yields
$\alpha(G)>10n/(11a)$.
Theorem~\ref{thm:counting} gives the stronger estimate
$\alpha(G)>n/a$ in terms of $a$.

\bigskip

\noindent\textbf{Notation}. 
We write $A\sqcup B$ for the disjoint union of two sets $A$ and $B$. For a positive integer $m$, let
$[m]=\{1,2,\ldots,m\}$.  For a vertex set $I\subseteq V(G)$, let
$N_G(I)=\{v\in V(G)\setminus I: N_G(v)\cap I\neq\emptyset\}$
and $N_G[I]=I\cup N_G(I)$. Two vertex sets $A,B\subseteq V(G)$ are \emph{anticomplete}
if there are no edges between $A$ and $B$.
We say that $A$ \emph{totally dominates} $B$ if every vertex of $B$
has a neighbor in $A$.
 For a vertex $v\in V(G)$ and a set $A\subseteq V(G)$,
let $N_A(v) = N_G(v)\cap A$ and $d_A(v)=|N_G(v)\cap A|$.

\smallskip

In Section 2, we collect the preliminary results regarding the augmentation of independent sets. In Section 3, we present the proof of Theorem \ref{thm:counting}. The paper concludes with several remarks in Section 4.

\section{Preliminaries}

An \emph{independent transversal} of disjoint vertex sets chooses exactly one
vertex from each set, with all chosen vertices being independent.
Independent transversals are useful for constructing independent sets
subject to prescribed constraints. King~\cite{King2011} used lopsided
independent transversals to prove that if  a graph $G$ satisfies
$\omega(G)>\frac{2}{3}(\Delta(G)+1)$, then $G$ contains an independent set
meeting every maximum clique of $G$.

Our main idea is an augmenting argument on independent sets underlying a lemma of King. Given a partition of $G$ and any vertex of $G$, the lemma describes
what happens when the independent transversals of the partition cannot
include a prescribed vertex: this failure yields a small set totally
dominating the prescribed vertex together with a union of some parts of the partition. Together with a simple counting estimate on the size of the neighbourhood of the small set, we  obtain Lemmas \ref{lem:single-root} and \ref{lem:pair-transversal}, which ensure that we can always perform the augmentation of the independent set.  
The following lemma is adapted from  the construction in
the proof of King~\cite[Lemma~8]{King2011}.

\begin{lem}[\cite{King2011}]\label{lem:augmentation}
Let $G$ be a graph with $V(G)=\{v\}\sqcup\left(\bigsqcup_{i\in[m]}V_i\right)$. Suppose that $V_1,\dots,V_m$ admit an independent transversal, but no such independent transversal is anticomplete to $v$. Then there exist a nonempty index set $S\subseteq[m]$ and disjoint independent sets $X,Y\subseteq V(G)$ with the following properties:
\begin{enumerate}[$(i)$, leftmargin = *]
\item $Y$ is an independent transversal of $V_i, i\in S$;
\item $v\in X$, and $X\setminus \{v\}\subseteq  \bigcup_{i\in S}V_i$ with $|X|\le |Y|$;
\item $\{v\}\cup\bigcup_{i\in S} V_i$ is totally dominated by $X\cup Y$ in $G$.
\end{enumerate}
\end{lem}

\smallskip

Given an independent set $I$, we define the \emph{edge redundancy} of $I$ as
\[
\eta(I) := \sum_{v\in N_G(I)}(d_I(v) - 1) = e_G(I,N_G(I))-|N_G(I)|. 
\]
It follows from the definition that at most $\eta(I)$ vertices of $G$ are adjacent to at least two vertices in $I$. Indeed, each such vertex contributes at least $1$ to the corresponding summand.

In the remainder of this section, we fix $a>0$ and assume that $G$ satisfies $\Delta(G)<6a/5-1$, together with an additional condition
\begin{equation}\label{eq:edge-bound}
\eta(I) < a|I|/5
 \qquad \text{for every nonempty independent set }I.
\end{equation}
 Specifically, for any two nonadjacent vertices $u,v\in V(G)$, we have
\begin{equation}\label{twocap}
|N_G(u) \cap N_G(v)|  = \eta(\{u,v\}) < \frac{2a}{5}.
\end{equation} 
\begin{remark}\label{induce:heredity}
The definition of $\eta$  implies that the bound \eqref{eq:edge-bound} is inherited by any  induced subgraph containing $I$ with the same constant
$a$. Indeed, the values $d_I(v)$ at the retained vertices  are unchanged  in the induced subgraph,
and only nonnegative summands are removed.
\end{remark}

Call a clique \emph{large} if it has more than $4a/5$ vertices.
Let $C_1,\ldots,C_m$ be pairwise disjoint maximal large cliques in $G$. For $b\in C_i$, put $N^\times(b)=N_G(b)\setminus C_i$. Since
$C_i\subseteq N_G[b]$,
\begin{equation}\label{eq:cross}
|N^\times(b)|=d_G(b)+1-|C_i|
 < \left( \frac{6}{5}a - 1\right) +1- \frac{4}{5}a = \frac{2}{5}a.
\end{equation}

Let $S\subseteq [m]$ , put
$U_S=\bigcup_{i\in S}C_i$ and $s=|S|$, and let
$T=\{b_i:i\in S\}$ be an independent transversal. Suppose $z$
is outside the cliques and anticomplete to $T$. Each vertex of
$U_S\setminus T$ receives one edge from $T$ within cliques. After these  edges are removed,
the number of the remaining edges between $U_S\setminus T$ and $T\cup \{z\}$ is
\begin{align}
|N_G(z)\cap U_S|
&+\sum_{b\in T}|N^\times(b)\cap U_S|\notag\\
&=e_G(T\cup\{z\},U_S\setminus T)-|U_S\setminus T|\notag\\
&\le e_G(T\cup\{z\},N_G(T\cup\{z\}))
      -|N_G(T\cup\{z\})|\notag\\
& = \eta(T\cup \{z\})\overset{\eqref{eq:edge-bound}}{<}\frac{a(s+1)}5.\label{eq:joint-count}
\end{align}
The first inequality follows from $U_S\setminus T\subseteq N_G(T\cup\{z\})$. By \eqref{eq:joint-count}, we also have
\begin{equation}\label{cliquetoout}
  |N_C(x)| < 2a/5  
\end{equation}
for every vertex $x$ outside a maximal clique $C$. 

\begin{lem}\label{lem:single-root}
For every vertex $z$ outside $\bigcup_{i=1}^m C_i$, the cliques
have an independent transversal anticomplete to $z$. 
\end{lem}

\begin{proof}
We proceed by induction on $m$. For $m=1$, the maximality of $C_1$ yields a vertex $v\in C_1$ that is nonadjacent to $z$. Thus $\{v\}$ already establishes the base case. Now suppose the result holds for $m-1$, where $m\ge2$.

Let $V_i=C_i\setminus N_G(z)$ for $i\le m-1$. By the maximality of $C_i$ for $i\le m-1$, we have $V_i\neq\emptyset$; moreover, there exists a vertex $v\in C_m$ that is not adjacent to $z$. Form an auxiliary graph $H$ on vertex set $\{v\}\cup\bigcup_{i\in[m-1]}V_i$ by retaining all original edges except those within each $V_i$. 
By induction, the cliques $C_i,1\le i\le m - 1$,  have an independent
transversal anticomplete to $z$, which is an independent
transversal of $V_i,1\le i\le m-1$, in $H$.
If $V_1,\ldots,V_{m-1}$ have an independent transversal $T$
in $H$ that is anticomplete to $v$, then $T\cup\{v\}$ is an
independent transversal of $C_1,\ldots,C_m$ anticomplete to $z$ in $G$.
If no such transversal $T$ exists, Lemma~\ref{lem:augmentation} gives $S,X,Y$
with $s=|S|\ge 1$. $Y$ is an independent transversal of $C_i,i\in S$, and is anticomplete to $z$. All vertices in $X$ belong to $U_{[m]}$, including $v\in C_m$, and $|X|\le s$.

Every vertex of the full union $U_S$ is either a neighbor of $z$
or a vertex totally dominated by $X\cup Y$ in $H$.
All auxiliary edges in $H$ involved are cross edges of $C_1,\dotsc,C_m$. Applying \eqref{eq:joint-count} with $Y$ and vertex $z$ and using
\eqref{eq:cross}, we obtain
\begin{align*}
|U_S|
&\le |N_G(z)\cap U_S|
 +\sum_{b\in Y}|N^\times(b)\cap U_S|
 +\sum_{u\in X}|N^\times(u)\cap U_S|\\
&<\frac{a(s+1)}5+\frac{2as}5
=\frac{a(3s+1)}5\le\frac{4as}5.
\end{align*}
This contradicts $|U_S|>4as/5$, since $|C_i|>4a/5$ for all $1\le i\le m$. Thus the induction is complete.
\end{proof}

\begin{lem}\label{lem:pair-transversal}
If $x,y$ are distinct nonadjacent vertices outside $\bigcup_{i=1}^m C_i$,
the cliques
have an independent transversal anticomplete to $\{x,y\}$. 
\end{lem}

\begin{proof}
Lemma~\ref{lem:single-root}, applied
to $x$, gives an independent transversal of the parts 
$V_i=C_i\setminus N_G(x),1\le i\le m$. Form the auxiliary graph on vertex set $\{y\}\sqcup \left(\bigsqcup_{i\in [m]} V_i\right)$ by deleting within-part edges and keeping all other
original edges. An independent transversal of $V_1,\dotsc,V_m$ anticomplete to $y$ in this graph is
the desired transversal in $G$. If no such independent transversal exists,  Lemma~\ref{lem:augmentation} gives $S,X,Y$. Put
$s=|S|\ge1$ and $X'=X\setminus\{y\}$.  $Y$ is an independent transversal of $C_i,i\in S$, and is anticomplete to $x$. Moreover,
$X'\subseteq U_{[m]}$ and $|X'|\le s-1$.

Every vertex of $U_S$ is covered by $N_G(x)$, a cross-neighborhood
of a member of $Y$, $N_G(y)$, or a cross-neighborhood of a member
of $X'$.  The first two
contributions together satisfy \eqref{eq:joint-count} with $Y$
and vertex $x$.
By Lemma~\ref{lem:single-root}, the cliques $C_i$, $i\in S$,
have an independent transversal anticomplete to $y$. Applying \eqref{eq:joint-count} to this transversal
gives $|N_G(y)\cap U_S|<a(s+1)/5$. The remaining contributions can be bounded by
\eqref{eq:cross}. Consequently
\begin{align*}
|U_S|
&\le |N_G(x)\cap U_S|
 +\sum_{b\in Y}|N^\times(b)\cap U_S|\\
&\quad +|N_G(y)\cap U_S|
 +\sum_{u\in X'}|N^\times(u)\cap U_S|\\
&<\frac{a(s+1)}5+\frac{a(s+1)}5+\frac{2a(s-1)}5
=\frac{4as}5.
\end{align*}
This contradicts $|U_S|>4as/5$.
\end{proof}

\section{Proof of Theorem~~\ref{thm:counting}}

\begin{proof}[Proof of Theorem~\ref{thm:counting}]
Fix $a>0$. If Theorem~\ref{thm:counting} fails, choose a counterexample
$G$ with the fewest vertices. Write  $k=\al(G)$, so
\begin{equation}\label{eq:counterexample}
n\ge ak,\qquad\om(G)<a,\qquad\Delta(G)<6a/5-1.
\end{equation}
For every nonempty proper vertex set $W$,
\begin{equation}\label{eq:hereditary}
|W|<a\al(G[W]).
\end{equation}
Indeed, $\om(G[W])\le\om(G)<a$ and
$\Delta(G[W])\le\Delta(G)<6a/5-1$, so $G[W]$ satisfies both
hypotheses with the same $a$. The minimality of $G$ therefore
gives \eqref{eq:hereditary}.

\bigskip
\noindent\textbf{Claim 1}. 
   $ \eta(I) < a|I|/5$
 for every nonempty independent set $I$.
\begin{proof}
    Let $I$ be an independent set of size $t$, and put
$R=V(G)\setminus N_G[I]$. Since $R$ is disjoint from and anticomplete to $I$, we have
$\alpha(G[R])+t\le k$. If $R$ is nonempty, it is proper, and
\eqref{eq:hereditary} yields
\[
|R|<a\al(G[R])\le a(k-t),\qquad
|N_G[I]|= n - |R|>at.
\]
If $R$ is empty, then $|N_G[I]| = n\ge ak\ge at$. Thus
$|N_G[I]|\ge at$ in every case. Consequently,
\begin{align*}
\eta(I) = e_G(I,N_G(I))-|N_G(I)|
&=\sum_{u\in I}(d_G(u)+1)-|N_G[I]|\notag\\
&\le t(\Delta(G)+1)-at\notag\\
&<\frac{6at}{5}-at=\frac{at}{5}.
\label{eq:derived-edge-bound}
\end{align*}
\end{proof}
This verifies condition \eqref{eq:edge-bound}. Hence, Lemmas~\ref{lem:single-root} and \ref{lem:pair-transversal} apply.

\bigskip
\noindent\textbf{Claim 2.} Distinct maximal large cliques in $G$ are disjoint. 

\begin{proof}
Suppose that two distinct maximal large cliques $C$ and $D$ intersect.
By maximality, there are nonadjacent vertices $u\in C\setminus D$
and $v\in D\setminus C$. Since $C\cap D\subseteq N_G(u)\cap N_G(v)$,
\eqref{twocap} gives $|C\cap D|<2a/5$. As $|C|,|D|>4a/5$, we have
\[
|C\cup D|=|C|+|D|-|C\cap D|>\frac{6a}{5}.
\]
Now choose $w\in C\cap D$. Then $C\cup D\subseteq N_G[w]$, so
\[
|C\cup D|\le \Delta(G)+1<\frac{6a}{5},
\]
a contradiction.
\end{proof}

Let $C_1,\ldots,C_h$ be distinct maximal large cliques and
$B_i\subseteq C_i$ cliques with $|B_i|>4a/5$. Let $Z$ consist
of one vertex, or two distinct nonadjacent vertices, all outside
$\bigcup_i C_i$. Consider
\[
H=G\left[\left(\bigcup_{i=1}^h B_i\right)\cup Z\right].
\]
Since $H$ is induced, $\om(H)\le\om(G)<a$ and
$\Delta(H)\le\Delta(G)<6a/5-1$. Moreover,
\eqref{eq:edge-bound} holds in $H$. For any $v\in V(H)\setminus B_i$, we have $v\notin C_i$.
By \eqref{cliquetoout}, $|N_G(v)\cap B_i|\le |N_G(v)\cap C_i|<2a/5<|B_i|$.
Thus $B_i$ is maximal in $H$. $B_1,\dotsc,B_h$ are disjoint,
so Lemma~\ref{lem:single-root} or
Lemma~\ref{lem:pair-transversal}, applied in $H$, gives an
independent transversal of the $B_i$ anticomplete to $Z$.
This remains true in $G$, since $H$ is induced. 

Maintain an independent set $I$, distinct maximal large cliques
$C_1,\ldots,C_h$, and large subcliques $B_i\subseteq C_i$ such that
\begin{equation}\label{eq:invariant}
|I|+h=k,\qquad I\cap B_i=\varnothing,
\qquad I\text{ is anticomplete to }B_i\quad(1\le i\le h).
\end{equation}
A maximum independent set $I$ initializes this configuration with
$h=0$. No vertex of $I$ belongs to a parent clique $C_i$. Indeed, it
would be distinct from every vertex of the nonempty $B_i$ and
adjacent to all of them. Thus
\[
I\subseteq R:=V(G)\setminus\bigcup_{i=1}^h C_i.
\]

Suppose $I\ne\varnothing$, and put $\ell=|I|=k-h$. Since
$C_1,C_2,\dotsc,C_h$ are pairwise disjoint and $\om(G)<a$,
\begin{equation}\label{eq:remainder}
|R|=n-\sum_{i=1}^h|C_i|\ge a(k-h)=a\ell.
\end{equation}
 We claim that $N_G[x]\cap I\ne\varnothing$ for every $x\in R$.
Otherwise choose $x\in R$ with $N_G[x]\cap I=\varnothing$.
It is outside $I$ and anticomplete to $I$. The preceding
subclique argument supplies an independent transversal $T$ of the $B_i$
anticomplete to $x$. Each $B_i$ is also disjoint from and
anticomplete to $I$, so $I\cup\{x\}\cup T$ is independent
of size $\ell+1+h=k+1 > \al(G)$, a contradiction.
For each $u\in I$, let 
\[
Q_u=\{x\in R:N_G[x]\cap I=\{u\}\}.
\]
We claim that $Q_u$ is a clique. Indeed,  if distinct $x,y\in Q_u$
were nonadjacent, neither could equal $u$, since every other
member of $Q_u$ is adjacent to $u$. Neither could be another
vertex of $I$, by the 
independence of $I$. Both are anticomplete to $I\setminus\{u\}$
and outside $\bigcup_iC_i$. The subclique argument gives
an independent transversal $T$ of all $B_i$ anticomplete to $x,y$. Then $(I\setminus\{u\})\cup\{x,y\}\cup T$
is independent of size $\ell-1+2+h=k+ 1 > \al(G)$, a contradiction.

 Let $M$ be the set of vertices in $R$ with at least two neighbors in $I$.
Since $N_G[x]\cap I\ne\varnothing$ for every $x\in R$,
the sets $M$ and $Q_u$ ($u\in I$) partition $R$.  By the definition of $\eta(I)$ and \eqref{eq:edge-bound}, 
\[
|M|\le \eta(I) <\frac{a\ell}{5}.
\]
Together with \eqref{eq:remainder}, this yields
\begin{equation}\label{eq:private-surplus}
\sum_{u\in I}|Q_u|=|R|-|M|>\frac{4a\ell}{5}.
\end{equation}
There are $\ell$ such cliques, so one, say $Q_u$, is large.
Extend it to a maximal large clique $C_{h+1}$. This differs from
every previous $C_i$, since it contains the nonempty set $Q_u$
lying entirely in $R$. Disjointness of distinct maximal large
cliques makes $C_{h+1}$ disjoint from $C_1,C_2,\dotsc,C_h$. Set $B_{h+1}=Q_u$ and replace $I$ by $I\setminus\{u\}$. By definition, $B_{h+1}$ is disjoint from and anticomplete to the new $I$. The sets $B_1,B_2,\dotsc,B_h$ retain these properties. Removing $u$ from $I$ and adding $C_{h+1}$ preserves $|I|+h=k$. Thus \eqref{eq:invariant} holds again. The step applies whenever $I$ is nonempty. After exactly $k$
steps it gives $I=\varnothing$ and $h=k$. At this point
\[
|R| = n-\sum_{i=1}^k|C_i|>n-ak\ge0.
\]
Choose a vertex $x\in R$. The same subclique argument supplies
an independent transversal of  $B_1,B_2,\dotsc,B_k$  anticomplete to $x$. Together
with $x$ it forms an independent set of size $k+1$, contradicting
$\al(G)=k$. This proves Theorem~\ref{thm:counting}. 
\end{proof}

\section{Concluding remarks}

Motivated by the chromatic profile studied by
Erd\H{o}s and Simonovits~\cite{ErdosSimonovits1973}
and explicitly formulated by Allen et al.~\cite{ABGKM},
we introduce an analogous parameter, which we call
the \emph{independence profile}. We define the \emph{independence profile}
by replacing colourability with the existence of a large
independent set. For a graph $H$ and a positive integer $\ell$, let
\[
\delta_\alpha(H,\ell)
:=\inf\left\{d\in[0,1]:
\begin{array}{l}
\alpha(G)\ge |V(G)|/\ell
\text{ for every }H\text{-free graph }G\\
\text{with }\delta(G)>d|V(G)|
\end{array}
\right\}.
\]
Since every $\ell$-colourable graph has an independent set of
size at least $|V(G)|/\ell$, we have
$
\delta_\alpha(H,\ell)\le\delta_\chi(H,\ell)$. Here $\delta_\chi(H,\ell)$ denotes the corresponding minimum-degree threshold for $\ell$-colourability. Our result, together with the sharpness constructions, determines
\[
\delta_\alpha(K_r,r-1)=\frac{5r-11}{5r-5}.
\]
For $r\ge4$, this gives a strict separation between the two profiles:
\[
\delta_\alpha(K_r,r-1)
=\frac{5r-11}{5r-5}
<
\frac{3r-7}{3r-4}
=\delta_\chi(K_r,r-1),
\]
where the last equality follows from the
Andr\'asfai--Erd\H{o}s--S\'os theorem.

It is natural to investigate the remaining values of the
independence profile. In particular, what are the exact values
of $\delta_\alpha(K_r,\ell)$ for $\ell>r-1$?
More generally, how does the independence profile of $H$-free
graphs compare with its chromatic profile?

\section*{Acknowledgments}

This work was supported by the National Key R\&D Program of China (2023YFA1010203), the
National Natural Science Foundation of China (12471336), and the Quantum Science and
Technology-National Science and Technology Major Project (2021ZD0302902).

\section*{Declaration on the use of AI}

During the preparation of this work, the authors used ChatGPT-6
to assist with mathematical exploration. The initial AI-generated
proof was difficult to follow and lacked references. The authors
therefore carefully examined the arguments, substantially revised
the proof for clarity and rigor, and supplied the relevant background
and references based on their own understanding of the problem.
The authors have reviewed and verified all proofs and arguments
and take full responsibility for the content of this paper.

\end{document}